\documentclass[12pt]{article}
\usepackage{amssymb,amsmath}
\usepackage{amsthm}
\usepackage{cases}
\usepackage{amsfonts}
\usepackage{cite,color,xcolor}
\usepackage[left=3.2cm,right=3.2cm,top=2.8cm,bottom=2.8cm]{geometry}
\usepackage[colorlinks,citecolor=blue,urlcolor=blue]{hyperref}
\usepackage[utf8]{inputenc}
 \usepackage[pagewise]{lineno}

\newtheorem{theorem}{Theorem}[section]
\newtheorem{corollary}{Corollary}[section]
\newtheorem{lemma}{Lemma}[section]

\newtheorem{remark}{Remark}[section]

\usepackage{txfonts}

\newcommand{\bal}{\begin{align}}
\newcommand{\bbal}{\begin{align*}}
\newcommand{\beq}{\begin{equation}}
\newcommand{\eeq}{\end{equation}}
\newcommand{\bca}{\begin{cases}}
\newcommand{\eca}{\end{cases}}

\newcommand{\pa}{\partial}
\newcommand{\fr}{\frac}

\newcommand{\De}{\Delta}

\newcommand{\cd}{\cdot}
\newcommand{\ep}{\varepsilon}
\newcommand{\dd}{\mathrm{d}}

\newcommand{\R}{\mathbb{R}}

\newcommand{\f}{\left}
\newcommand{\g}{\right}

\begin{document}
\title{Global existence and blow-up for the Euler-Poincar\'{e} equations with a class of initial data}

\author{
Jinlu Li$^{1}$\thanks{School of Mathematics and Computer Sciences, Gannan Normal University, Ganzhou 341000, China. E-mail: \texttt{lijinlu@gnnu.edu.cn}},
Yanghai Yu$^{2}$\thanks{School of Mathematics and Statistics, Anhui Normal University, Wuhu 241002, China.
E-mail: \texttt{yuyanghai214@sina.com}(Corresponding author)}
and Weipeng Zhu$^{3}$\thanks{School of Mathematics and Big Data, Foshan University, Foshan, Guangdong 528000, China.
E-mail: \texttt{mathzwp2010@163.com}}
}

\date{\today}

\maketitle\noindent{\hrulefill}

{\bf Abstract:} In this paper we  investigate the Cauchy problem of d-dimensional Euler-Poincar\'{e} equations. By choosing a class of new and special initial data, we can transform this d-dimensional Euler-Poincar\'{e} equations into the Camassa-Holm type equation in the real line. We first obtain some global existence results and then present a new blow-up result to the system under some different assumptions  on this special class of initial data.


\vskip0mm\noindent{\hrulefill}

\section{Introduction}

In this paper, we consider the Cauchy problem of the following higher dimensional Euler-Poincar\'{e} equations
\begin{equation}\label{e0}
\partial_tm+u\cdot \nabla m+(\nabla u)^{\top}\cd m+(\mathrm{div} u)m=0, \quad
m=(1-\De)u,
\end{equation}
where $t\in \R^+$ and $x=(x_1,x_2,\cdots,x_d)\in \R^d$. $u(t,x)=(u_1,u_2,\cdots,u_d)(t,x)$ denotes the velocity of the fluid and $m(t,x)=(m_1,m_2,\cdots,m_d)(t,x)$ represents the momentum. The notation $(\nabla u)^{\top}$
denotes the transpose of the matrix $\nabla u$. It is useful to recast
equation \eqref{e0} in the component-wise form as
$$\partial_tm_i+\sum^d_{j=1}u_j\partial_{x_j}m_i+\sum^d_{j=1}(\partial_{x_i}u_j)m_j+m_i\sum^d_{j=1}\partial_{x_j}u_j=0.$$
While the Euler-Poincar\'{e} equations \eqref{e0} is also called the higher dimensional Camassa-Holm equations. In one dimension, the system \eqref{e0} is the classical Camassa-Holm (CH) equation
$$\partial_tm+u\pa_x m+2m\pa_xu=0, \quad
m=(1-\pa_{xx})u.
$$

Like the KdV equation, the CH equation  describes the unidirectional propagation of waves at the free surface of shallow water under the influence of gravity \cite{Camassa,Camassa.Hyman,Constantin.Lannes}.  It is completely integrable \cite{Camassa,Constantin-P}, has a bi-Hamiltonian structure \cite{Constantin-E,Fokas}, and admits exact peaked solitons of the form $ce^{-|x-ct|}$($c>0$), which are orbitally stable \cite{Constantin.Strauss}. It is worth mentioning that the peaked solitons present the characteristic for the travelling water waves of greatest height and largest amplitude and arise as solutions to the free-boundary problem for incompressible Euler equations over a flat bed, cf. \cite{Constantin2,Constantin.Escher4,Constantin.Escher5,Toland}. The local well-posedness and ill-posedness for the Cauchy problem of the CH equation in Sobolev spaces and Besov spaces was discussed in \cite{Constantin.Escher,Constantin.Escher2,d1,d3,Guo-Yin,Li-Yin,Ro}. It was shown that there exist global strong solutions to the CH equation \cite{Constantin,Constantin.Escher,Constantin.Escher2} and finite time blow-up strong solutions to the CH equation \cite{Constantin,Constantin.Escher,Constantin.Escher2,Constantin.Escher3}.  The global conservative and dissipative solutions of CH equation were discussed in \cite{Bressan.Constantin,Bressan.Constantin2}.

The Euler-Poincar\'{e} equations arise in diverse scientific applications and enjoy several remarkable properties
both in the one-dimensional and multi-dimensional cases. The Euler-Poincar\'{e} equations were first studied by Holm, Marsden, and Ratiu in
1998 as a framework for modeling and analyzing fluid dynamics \cite{H-M-R1,H-M-R2}, particularly
for nonlinear shallow water waves, geophysical fluids and turbulence modeling.
Later, the Euler-Poincar\'{e} equations have many further interpretations beyond fluid applications. For instance, in 2-D, it is exactly the same as the averaged template matching equation for computer vision \cite{H-M-A}. Also, the Euler-Poincar\'{e} equations  have important applications in computational anatomy, it can be regarded as an evolutionary equation for a geodesic motion on a diffeomorphism group and it is associated with Euler-Poincar\'{e} reduction via symmetry (see, e.g, \cite{H-S-C,Y}).

 The rigorous analysis of the Euler-Poincar\'{e} equations  with $d\geq1$ was initiated by Chae and Liu \cite{Chae.Liu} who obtained the local well-posedness in Hilbert spaces $m_0\in H^{s+\frac d2},\ s\geq2$ and also gave a  blow-up criterion, zero $\alpha$ limit and the Liouville type theorem. Li, Yu and Zhai \cite{L.Y.Z} proved that the solution to higher dimensional Camassa-Holm equations with a large class of smooth initial data blows up in finite time or exists globally in time, which reveals the nonlinear depletion mechanism hidden in the Euler-Poincar\'{e} system. Luo and Yin \cite{Luo-Yin} obtained a new blow-up result to the periodic Euler-Poincar\'{e} system for a special class of smooth initial data by using the rotational invariant properties of the system.
By means of the Littlewood-Paley theory, Yan and Yin \cite{Y.Y} established the local existence and uniqueness in Besov spaces $B^s_{p,r}$ with $s>\max\{\frac32,1+\frac dp\}$ and $s=1+\frac dp,\ 1\leq p\leq 2d,\ r=1$.
For more results of the Euler-Poincar\'{e} equations, we refer the reads to see \cite{L.D.Z,LDL,LiG,Z.Y.L}.
In this paper, we consider the Euler-Poincar\'{e} equations in $\R^d$ with a class of special initial data be of the form $f(x_1+x_2+\cdots x_d)$ and study the global existence and blow-up property of the corresponding solution under certain assumptions on this initial data.

\subsection{Reduction of System}\label{subsec1.1}
Let us assume that
\begin{equation}\label{1d}
\begin{cases}
\pa_tn+dv\pa_xn+2dv_xn=0, \quad (t,x)\in \R^+\times \R,\\
n=(1-d\pa^2_x)v,\\
n(0,x)=n_0(x),\quad x\in  \R.
\end{cases}
\end{equation}
Then we can find that
\begin{align*}
&m(t,x)=n(t,x_1+x_2+\cdots +x_d)\vec{e}\quad\text{with}\quad\vec{e}:=\big(\underbrace{1,\cdots,1}_{d}\big),\\
&u(t,x)=v(t,x_1+x_2+\cdots +x_d)\vec{e},
\end{align*}
satisfy the following higher dimensional Euler-Poincar\'{e} equations
\begin{equation}\tag{EP}\label{EP}
\begin{cases}
\partial_tm+u\cdot \nabla m+(\nabla u)^{\top}\cd m+(\mathrm{div} u)m=0, \quad& (t,x)\in \R^+\times \R^d,\\
m=(1-\De)u,\\
u(0,x)=u_0(x),\quad& x\in \R^d.
\end{cases}
\end{equation}
Precisely speaking, the Cauchy Problem of higher dimensional Euler-Poincar\'{e} equations  \eqref{EP} is transformed into that of the new Camassa-Holm type equation \eqref{1d} in the real-line. Furthermore, we can rewrite the equation  \eqref{1d} as follows
\begin{equation}\tag{d-CH}\label{d-CH}
\begin{cases}
\pa_tv+dvv_x=-\pa_x(1-d\pa^2_x)^{-1}\f(dv^2+\frac{d^2}{2}(\pa_xv)^2\g) , \quad& (t,x)\in \R^+\times \R,\\
v(0,x)=v_0(x),\quad& x\in \R.
\end{cases}
\end{equation}
From now, we mainly focus on the equation  \eqref{d-CH}.
\subsection{Main results}\label{subsec1.2}

First, following the  similar proofs for the Camassa-Holm equation in \cite{Ro,d1}, we can obtain the well-posedness result for \eqref{d-CH}  as follows.

\begin{theorem}[Local well-posedness]\label{th0}
 Let $v_0 \in H^s(\mathbb{R})$ with $s>\frac{3}{2}$. Then there exists a time $T^*>0$ and a unique strong solution $v(t,x) \in C\left([0, T^*) ; H^s(\mathbb{R})\right) \cap C^1\left([0, T^*) ; H^{s-1}(\mathbb{R})\right)$ of the Cauchy problem  \eqref{d-CH}.
\end{theorem}

\begin{remark}
Let $v_0 \in H^s(\mathbb{R})$ with $s>2$. Assume that $u_0(x)=v_0(x_1+x_2+\cdots +x_d)\vec{e}$, by the embedding 
$
H^s(\mathbb{R})\hookrightarrow \mathcal{C}^{s-\fr12-\ep}(\mathbb{R}),
$
where we take $\ep=0$ if $s-\fr12\notin \mathbf{Z}^+$ and $\ep=0^+$ if $s-\fr12\in \mathbf{Z}^+$, then we have $u_0(x)\in \mathcal{C}^{s-\fr12-\ep}(\R^d)$. Moreover, we can obtain the unique solution $u(t,x)=v(t,x_1+x_2+\cdots +x_d)\vec{e}\in C([0,T^*);\mathcal{C}^{s-\fr12-\ep}(\R^d))$ of the Cauchy problem \eqref{EP} with initial data $u_0(x)$. Based on the above observation, we can transform the global existence or blow-up property of solutions for the Cauchy problem \eqref{EP}  into that for the Cauchy problem \eqref{d-CH}.
\end{remark}

Finally, we state the following two global existence results and one blow-up result.

\begin{theorem}[Global existence]\label{th1}
Let $n_0=v_0-d\pa^2_xv_0$ with $d\geq 1$. Assume that $v_0\in H^s(\mathbb{R})$ with $s>\frac{3}{2}$ satisfying one of the following
\begin{description}
  \item[(1)] $n_0\geq 0$ on $(-\infty,+\infty)$;
  \item[(2)] $n_0\leq 0$ on $(-\infty,x_0]$ and $n_0\geq 0$ on $[x_0,+\infty)$ for some point $x_0\in\R$.
\end{description}
Then the corresponding solution of the Cauchy problem \eqref{d-CH} exists globally in time.
\end{theorem}
\begin{corollary}\label{cor010}
Let $d\geq 1$. Assume that $u_0(x_0)=v_0(x_1^0+x_2^0+\cdots+x_d^0)\vec{e}$ with $v_0$ satisfying the assumption of Theorems \ref{th1}. Then the corresponding solution of the Euler-Poincar\'{e} equations \eqref{EP} exists globally in time.
\end{corollary}
\begin{theorem}[Formation of singularities]\label{th3}
Let $d\geq 1$. Assume that $v_0\in H^s(\mathbb{R})$ with $s>\frac{3}{2}$ satisfying that for some point $x_0\in \R$
\bal\label{as}
v'_0(x_0)<-\frac{1}{\sqrt{d}}|v_0(x_0)|.
\end{align}
Then the corresponding solution of the Cauchy problem \eqref{d-CH} blows up in finite time.
\end{theorem}
\begin{corollary}\label{cor10}
Let $d\geq 1$. Assume that $u_0(x_0)=v_0(x_1^0+x_2^0+\cdots+x_d^0)\vec{e}$ with $v_0$ satisfying the assumption of Theorems \ref{th3}. Then the corresponding solution of the Euler-Poincar\'{e} equations \eqref{EP} blows up in finite time.
\end{corollary}
When assuming that the initial condition is stronger than \eqref{as}, we have
\begin{corollary}\label{cor1}
Let $d\geq 1$. Assume that $v_0\in H^s(\mathbb{R})$ with $s>\frac{3}{2}$ satisfying $n_0=v_0-d\pa^2_xv_0\geq 0$ on $(-\infty,x_0]$ and $n_0\leq 0$ on $[x_0,+\infty)$ for some point $x_0\in\R$ and $n_0$ changes sign. Then the corresponding solution of the Cauchy problem \eqref{d-CH} blows up in finite time.
\end{corollary}
\section{Preliminaries}\label{sec2}

Now let us consider the associated Lagrangian flow
\begin{equation}\label{L}
\begin{cases}
\frac{\dd q}{\dd t}=dv(t,q(t,x)),\quad (t,x)\in[0, T )\times\R,\\
q(0,x)=x,~~~~~~x\in\mathbb{R}.
\end{cases}
\end{equation}
Following the similar proofs for Camassa-Holm equation in \cite{Constantin}, we can obtain the following results:
\begin{lemma}\label{lm1} Assume that $v_0 \in H^s(\R)$ with $s>\frac{3}{2}$ and $T$ is the existence time of the corresponding solution of \eqref{d-CH}. Then one has, with $n=(1-d\pa^2_x)v$,
$$
n(t, q(t, x))q_x^2(t, x)=n_0(x), \quad (t,x)\in[0, T )\times\R,
$$
where $q(t, x) \in C^1([0, T), \R)$ is a unique solution of the Cauchy problem \eqref{L} and satisfies that
$$
q_x(t, x)=\exp \left(d\int_0^t v_x(\tau, q(\tau, x)) \dd \tau\right)>0,\quad t \in[0, T) .
$$
\end{lemma}

\begin{lemma}\label{lm2} Let $v_0 \in H^s(\mathbb{R})$ with $s>\frac{3}{2}$ and $T$ is the maximal existence time of the corresponding solution of \eqref{d-CH}. Then the corresponding solution of \eqref{d-CH} blows up in finite time if and only if
$$
\lim _{t \rightarrow T} \inf \left(\inf _{x \in \mathbb{R}} v_x(t, x)\right)=-\infty .
$$
That is, singularities can arise only in the form of wave breaking.
\end{lemma}

\begin{lemma}\label{lm3} Let $v_0 \in H^s(\mathbb{R})$ with $s>\frac{3}{2}$ and $T$ is the maximal existence time of the corresponding solution of \eqref{d-CH}. Then for any $t \in[0, T)$, we have
\bbal
\|v(t)\|^2_{L^2}+d\|\pa_xv(t)\|^2_{L^2}=\|v_0\|^2_{L^2}+d\|\pa_xv_0\|^2_{L^2}\leq 2d\|v_0\|^2_{H^1}.
\end{align*}
\end{lemma}

\section{Proof of Theorems}\label{sec3}
In this section, inspired by \cite{Constantin,Brandolese}, we shall prove our main results. The proof of Theorem \ref{th0} is standard, we omit the details. Next, we begin to prove the remaining Theorems.
\subsection{Proof of Theorem \ref{th1}}\label{subsec3.1}
\quad{\bf Case 1:} $n_0\geq 0$ on $(-\infty,+\infty)$. \\
By the assumption of $n_0$ and Lemma \ref{lm1}, one has
$$n(t,x)\geq 0,\quad (t,x)\in[0, T )\times\R.$$
Since $
\left(1-d \partial_x^2\right)^{-1} f=p_d \ast f$ for any $f \in L^2(\R)$ with $p_d(x)=\frac{1}{2 \sqrt{d}} \exp \left(-\frac{|x|}{\sqrt{d}}\right)$, one has
\bbal
v(t,x)=\frac{1}{2 \sqrt{d}} e^{-\frac{|x|}{\sqrt{d}}}*n(t,x)=\frac{1}{2 \sqrt{d}} \int^{\infty}_{-\infty}e^{-\frac{|x-\xi|}{\sqrt{d}}}n(t,\xi)\dd\xi.
\end{align*}
Thus due to Lemma \ref{lm3}, we have for $(t,x)\in[0, T) \times \R$
\bbal
|\pa_xv(t,x)|&=\f|\frac{1}{2 d} \int^{\infty}_{-\infty}\mathrm{sgn}(x-\xi)e^{-\frac{|x-\xi|}{\sqrt{d}}}n(t,\xi)\dd\xi\g|\\
&\leq \frac{1}{2 d} \int^{\infty}_{-\infty}e^{-\frac{|x-\xi|}{\sqrt{d}}}n(t,\xi)\dd\xi= \frac{1}{\sqrt{d}} v(t,x)\\
&\leq \sqrt{\frac{\|v(t,x)\|_{L^2}^2+d\|\pa_x v(t,x)\|_{L^2}^2}{2d}} \leq \|v_0\|_{H^1}.
\end{align*}
In view of Lemma \ref{lm2}, this shows the existence time $T=\infty$ and completes the proof of
Theorem \ref{th1}.

{\bf Case 2:} $n_0\leq 0$ on $(-\infty,x_0]$ and $n_0\geq 0$ on $[x_0,+\infty)$ for some point $x_0\in\R$.\\
Since $q(t,x)$ is an increasing diffeomorphism of $\R$ for $t\in[0, T )$, we deduce from Lemma \ref{lm1} that
\begin{equation}
\begin{cases}
n(t,x)\leq 0,  \quad \text{if}\; \ x\leq q(t,x_0),\\
n(t,x)\geq 0, \quad\text{if}\;  \ x\geq q(t,x_0),
\end{cases}
\end{equation}
and $n(t,q(t,x_0))=0$.
Then, we have for $x\geq q(t,x_0)$
\bal\label{hyl1}
v_x(t,x)=-\frac{1}{\sqrt{d}}v(t,x)+\frac{1}{d}e^{\frac{x}{\sqrt{d}}}\int^{+\infty}_xe^{-\frac{\xi}{\sqrt{d}}}n(t,\xi)\dd \xi\geq -\frac{1}{\sqrt{d}}v(t,x),
\end{align}
while $x\leq q(t,x_0)$
\bal\label{hyl2}
v_x(t,x)=\frac{1}{\sqrt{d}}v(t,x)-\frac{1}{d}e^{\frac{-x}{\sqrt{d}}}\int_{-\infty}^xe^{-\frac{\xi}{\sqrt{d}}}n(t,\xi)\dd \xi\geq \frac{1}{\sqrt{d}}v(t,x).
\end{align}
It thus follows from the relations \eqref{hyl1}-\eqref{hyl2} and Lemma \ref{lm3} that
\bbal
v_x(t,x)\geq -\frac{1}{\sqrt{d}}|v(t,x)|\geq -\|v_0\|_{H^1},\quad (t,x)\in[0, T) \times \R.
\end{align*}
In view of Lemma \ref{lm2}, this shows the existence time $T=\infty$ and completes the proof of
Theorem \ref{th3}.

\subsection{Proof of Theorem \ref{th3}}\label{subsec3.3}

We only need to show that blow-up results hold for initial data $v_0\in H^3(\R)$. Then the continuous dependence on initial data ensures the validity for all $H^s(\mathbb{R})$ with $s>\frac{3}{2}$. Let $v \in C\left([0, T), H^3(\R)\right) \cap C^1\left([0, T), H^2(\R)\right)$ be the solution of \eqref{d-CH}. Since $p_d(x)=\frac{1}{2 \sqrt{d}} \exp \left(-\frac{|x|}{\sqrt{d}}\right)$, one has
$$
\left(1-d \partial_x^2\right)^{-1} f=p_d \ast f, \quad f \in L^2(\R) .
$$
Differentiating Eq $\eqref{d-CH}_1$ with respect to $x$ yields
\bal\label{v1}
\pa_tv_x+dvv_{xx}=-\frac d2 v^2_x+v^2-p_d \ast\f(v^2+\frac d2(\pa_xv)^2\g).
\end{align}
We introduce two new functions $$w:=\frac{1}{\sqrt{d}}v\quad\text{and}\quad V:=v_x.$$
Then we can obtain from \eqref{v1} that
\bbal
&\pa_tV+dvV_{x}=-\frac d2 V^2+dw^2-dp_d \ast\f(w^2+\frac 12V^2\g),
\\&\pa_tw+dvw_{x}=-d^{\fr32}\pa_xp_d \ast\f(w^2+\frac 12V^2\g).
\end{align*}
From \eqref{L} and the above, we deduce  that
\bbal
&\frac{\dd}{\dd t}V(t,q(t,x))=\f(-\frac d2 V^2+dw^2-dp_d \ast\f(w^2+\frac 12V^2\g)\g)(t,q(t,x)),\\
&\frac{\dd}{\dd t}w(t,q(t,x))=-d^{\fr32}\f(\pa_xp_d \ast\f(w^2+\frac 12V^2\g)\g)(t,q(t,x)),
\end{align*}
which implies directly that
\bbal
&\frac{\dd}{\dd t}\f(w+V\g)(t,q(t,x))=d\f(-\frac 12 V^2+w^2-F(w,V)\g)(t,q(t,x)),\\
&\frac{\dd}{\dd t}\f(w-V\g)(t,q(t,x))=-d\f(-\frac 12 V^2+w^2-G(w,V)\g)(t,q(t,x)),
\end{align*}
where we set
\bbal
&F(w,V):=p_d*\f(w^2+\frac 12V^2\g)+\sqrt{d}\pa_xp_d*\f(w^2+\frac 12V^2\g),\\
&G(w,V):=p_d*\f(w^2+\frac 12V^2\g)-\sqrt{d}\pa_xp_d*\f(w^2+\frac 12V^2\g).
\end{align*}
Now we claim that
\bal\label{0}
F(w,V)(t,x)\geq \frac12w^2(t,x)\quad\text{and}\quad G(w,V)(t,x)\geq \frac12w^2(t,x),
\end{align}
which in turn gives that
\bal\label{00}
p_d*\f(w^2+\frac 12V^2\g)=\fr12\f(F(w,V)(t,x)+G(w,V)\g)\geq \frac12w^2(t,x).
\end{align}
In fact, one has
\bal\label{p0}
p_d*\f(w^2+\frac 12V^2\g)&=\frac{1}{2\sqrt{d}}e^{\frac{x}{\sqrt{d}}}\int^{+\infty}_xe^{-\frac{\xi}{\sqrt{d}}}\f(w^2+\frac 12V^2\g)\dd \xi\nonumber\\
&\quad+\frac{1}{2\sqrt{d}}e^{\frac{-x}{\sqrt{d}}}\int_{-\infty}^xe^{\frac{\xi}{\sqrt{d}}}\f(w^2+\frac 12V^2\g)\dd \xi
\end{align}
and
\bal\label{p1}
\sqrt{d}\pa_xp_d*\f(w^2+\frac 12V^2\g)&=\frac{1}{2\sqrt{d}}e^{\frac{x}{\sqrt{d}}}\int^{+\infty}_xe^{-\frac{\xi}{\sqrt{d}}}\f(w^2+\frac 12V^2\g)\dd \xi\nonumber\\
&\quad-\frac{1}{2\sqrt{d}}e^{\frac{-x}{\sqrt{d}}}\int_{-\infty}^xe^{\frac{\xi}{\sqrt{d}}}\f(w^2+\frac 12V^2\g)\dd \xi.
\end{align}
On the one hand, performing \eqref{p0}-\eqref{p1} yields that
\bbal
F(w,V)&=\frac{1}{\sqrt{d}}e^{\frac{x}{\sqrt{d}}}\int^{+\infty}_xe^{-\frac{\xi}{\sqrt{d}}}\f(w^2+\frac 12V^2\g)\dd \xi.
\end{align*}
It is not difficult to verify that
\bbal
\frac{1}{\sqrt{d}}e^{\frac{x}{\sqrt{d}}}\int^{+\infty}_xe^{-\frac{\xi}{\sqrt{d}}}\f(w^2+ V^2\g)\dd \xi
&\geq -\frac{1}{\sqrt{d}}e^{\frac{x}{\sqrt{d}}}\int^{+\infty}_xe^{-\frac{\xi}{\sqrt{d}}}\f(2w V\g)\dd \xi
\\&\geq w^2-\frac{1}{\sqrt{d}}e^{\frac{x}{\sqrt{d}}}\int^{+\infty}_xe^{-\frac{\xi}{\sqrt{d}}}w^2\dd \xi,
\end{align*}
which implies
\bal\label{1}
F(w,V)\geq \frac12w^2.
\end{align}
On the other hand, performing \eqref{p0}--\eqref{p1} yields that
\bbal
G(w,V)&=\frac{1}{\sqrt{d}}e^{\frac{-x}{\sqrt{d}}}\int^{x}_{-\infty}e^{\frac{\xi}{\sqrt{d}}}\f(w^2+ \fr12V^2\g)\dd \xi.
\end{align*}
Similarly,
\bbal
\frac{1}{\sqrt{d}}e^{\frac{-x}{\sqrt{d}}}\int^{x}_{-\infty}e^{\frac{\xi}{\sqrt{d}}}\f(w^2+V^2\g)\dd \xi
&\geq \frac{1}{\sqrt{d}}e^{\frac{-x}{\sqrt{d}}}\int^{x}_{-\infty}e^{\frac{\xi}{\sqrt{d}}}\f(2w V\g)\dd \xi
\\&\geq w^2-\frac{1}{\sqrt{d}}e^{-\frac{x}{\sqrt{d}}}\int^{x}_{-\infty}e^{\frac{\xi}{\sqrt{d}}}w^2\dd \xi,
\end{align*}
which implies
\bal\label{2}
\frac{1}{\sqrt{d}}e^{\frac{-x}{\sqrt{d}}}\int_{-\infty}^xe^{\frac{\xi}{\sqrt{d}}}\f(w^2+\frac 12V^2\g)\dd \xi\geq \frac12w^2.
\end{align}
Combing \eqref{1} and \eqref{2} yields \eqref{0}.

Using \eqref{00}, we have
\bal\label{100}
&\frac{\dd}{\dd t}V(t,q(t,x))\leq \f(-\frac d2 V^2+dw^2-\frac12dw^2\g)(t,q(t,x))\leq \frac d2\f(w^2-V^2\g)(t,q(t,x)),
\end{align}
and using \eqref{0}
\bal
&\frac{\dd}{\dd t}\f(w+V\g((t,q(t,x))\leq \frac d2\f(w^2-V^2\g)(t,q(t,x)),\label{101}\\
&\frac{\dd}{\dd t}\f(w-V\g((t,q(t,x))\geq -\frac d2\f(w^2-V^2\g)(t,q(t,x)).\label{102}
\end{align}
Solving the above inequalities,  we obtain
\bbal
&\frac{\dd}{\dd t}\Big(e^{\frac d2\int^t_0B(\tau,x)\dd \tau}A(t,x)\Big)\geq 0,\\
&\frac{\dd}{\dd t}\Big(e^{-\frac d2\int^t_0A(\tau,x)\dd \tau}B(t,x)\Big)\leq 0,
\end{align*}
where, for simplicity we set
\bbal
A(t,x)=(w-V)(t,q(t,x)) \quad\text{and}\quad B(t,x)=(w+V)(t,q(t,x)).
\end{align*}
Since $B(0,x_0)<0<A(0,x_0)$,  then we have $B(t,x_0)<0<A(t,x_0)$. Moreover, we obtain from \eqref{101}-\eqref{102} that
\bbal
&A(t,x_0)\geq e^{\frac d2\int^t_0B(\tau,x_0)\dd \tau}A(t,x_0) \geq A(0,x_0)>0, \\
&B(t,x_0)\leq e^{-\frac d2\int^t_0A(\tau,x_0)\dd \tau}B(t,x_0)\leq B(0,x_0)<0.
\end{align*}
For simplicity, we denote $$g(t):=V(t,q(t,x_0)),\quad A(t):=A(t,x_0)\quad\text{and}\quad B(t):=B(t,x_0).$$
Then, from \eqref{100}, we have
\bal\label{g}
g'(t)\leq \frac d2 A(t)B(t)\leq \frac d2 A(0)B(0).
\end{align}

Now we assume that the solution $v(t)$ of \eqref{1d} exists globally in time $t \in[0, \infty)$, that is, $T=\infty$. We next show  this leads to a contradiction.

Integrating \eqref{g} on $[0, t)$ yields
\bal\label{g1}
g(t) \leq g(0)+\frac d2 A(0)B(0)t .
\end{align}
Since $A(0) B(0)<0$, we have
$$
\lim _{t \rightarrow \infty} g(t)=-\infty .
$$
But
$$
\|v(t)\|_{\infty} \leq \|v_0\|_{H^1}<\infty, \quad t \in[0, \infty) .
$$
Hence there exists some $t_0>0$ such that
$$
g^2(t) \geq \frac{2}{d}\|v_0\|_{H^1}^2, \quad t \in\left[t_0, \infty\right).
$$
Combining the latter inequality with \eqref{g}, we have derived the inequality
\bbal
g'(t)\leq \frac d2 A(t)B(t)=-\frac d2 g^2(t)+\frac 12 v^2 \leq -\frac d4g^2(t), \quad t \in\left[t_0, \infty\right).
\end{align*}
On the other hand, by the assumptions \eqref{as}, we have
$$
g(0)=v'_0(x_0)<0.
$$
It then follows from \eqref{g1} that $g(t)<0$ for all $t \geq 0$. The differential inequality \eqref{g} can be therefore solved easily for the solution $g(t), t \in\left[t_0, \infty\right)$, that is
$$
\frac{1}{g\left(t_0\right)}-\frac{1}{g(t)}+\frac{d}{4}\left(t-t_0\right) \leq 0, \quad t \in\left[t_0, \infty\right).
$$
Since $-\frac{1}{g(t)}>0$, we have
$$
\frac{1}{g\left(t_0\right)}+\frac{d}{4}\left(t-t_0\right)<\frac{1}{g\left(t_0\right)}-\frac{1}{g(t)}+\frac{d}{4}\left(t-t_0\right)<0, \quad t \in\left[t_0, \infty\right),
$$
which leads to a contradiction as $t \rightarrow \infty$. This shows that $T<\infty$ and the proof of Theorem \ref{th3} is complete.

\section*{Acknowledgments}
J. Li is supported by the National Natural Science Foundation of China (12161004), Training Program for Academic and Technical Leaders of Major Disciplines in Ganpo Juncai Support Program(20232BCJ23009) and Jiangxi Provincial Natural Science Foundation (20224BAB201008). Y. Yu is supported by the National Natural Science Foundation of China (12101011). W. Zhu is supported by the National Natural Science Foundation of China (12201118) and Guangdong Basic and Applied Basic Research Foundation (2021A1515111018).

\section*{Declarations}
\noindent\textbf{Data Availability} No data was used for the research described in the article.

\noindent\textbf{Conflict of interest}
The authors declare that they have no conflict of interest.

\end{document}